\documentclass[a4paper,11pt]{article}
\RequirePackage[T1]{fontenc}
\usepackage{palatino}
\usepackage{hyperref}
\hypersetup{colorlinks=true, linkcolor=blue,  anchorcolor=blue, citecolor=blue, filecolor=blue, menucolor=blue, urlcolor=blue}
\usepackage[english]{babel}
\usepackage[utf8]{inputenc}
\usepackage[font={footnotesize,sf},labelfont=bf]{caption}
\usepackage{amsmath,amssymb,amsfonts,mathrsfs,amsthm}
\usepackage{graphicx,pdfpages,enumitem,cite}
\theoremstyle{plain}% default 

\theoremstyle{plain} 
 
\theoremstyle{definition} 
 
\usepackage{geometry}
\title{Data-driven nonlinear model reduction to spectral submanifolds in mechanical systems}
\author{M. Cenedese$^{1}$, J. Axås$^{1}$, H. Yang$^{2}$, M. Eriten$^{2}$ and G. Haller$^{1}$
\vspace{3.5mm}\\  %\thanks{Corresponding author: \href{mailto:georgehaller@ethz.ch}{georgehaller@ethz.ch}}
$^{1}$Institute for Mechanical Systems, ETH Z\"urich\\ Leonhardstrasse 21, 8092 Z\"urich, Switzerland \vspace{1.5mm}\\
$^{2}$Department of Mechanical Engineering, University of Wisconsin - Madison\\ 1513 University Avenue, Madison, WI 53706, USA}
\date{\today}
\begin{document}
\maketitle
\begin{abstract}
	\noindent While data-driven model reduction techniques are well-established for linearizable mechanical systems, general approaches to reducing non-linearizable systems with multiple coexisting steady states have been unavailable. In this paper, we review such a data-driven nonlinear model reduction methodology based on spectral submanifolds (SSMs). As input, this approach takes observations of unforced nonlinear oscillations to construct normal forms of the dynamics reduced to very low dimensional invariant manifolds. These normal forms capture amplitude-dependent properties and are accurate enough to provide predictions for non-linearizable system response under the additions of external forcing. We illustrate these results on examples from structural vibrations, featuring both synthetic and experimental data.
\end{abstract}

\section{Introduction}
Dimensionality reduction for data sets representing high-dimensional nonlinear mechanical systems is of crucial importance in science and technology. Low-dimensional models are expected to reduce computational cost and capture the essential physics of a high-dimensional system from data. Given the growing interest, for example, in light mechanical structures and MEMS devices, there is the need for truly nonlinear models, capturing amplitude-dependent properties and com-peting steady states solutions, which are increasingly important to identify, as highlighted in experiments of nonlinear mechanical vibrations \cite{Peeters2011a,Peeters2011b,Kurt2014,Ehrhardt2016,Renson2016b,Moore2019,TRC2019_PartI,TRC2019_PartII}. Predicting coexisting stable and unstable forced responses for a broad range of forcing amplitudes and frequencies is paramount in structural dynamics. However, a generally applicable technique returning such reliable low-dimensional models of nonlinear mechanical vibrations has not emerged yet.

The most common approaches to data-driven reduced-order modeling are the Proper Orthogonal Decomposition (POD) followed by a Galerkin projection \cite{Awrejcewicz2004,Holmes2012,Lu2019,Hijazi2020} and the Dynamic Mode Decomposition (DMD) \cite{Schmid2010,Kutz2016}. The former approach requires the knowledge of the governing equations of motion and, once a relevant number of modes is identified from data, projects these equations onto those modes to construct a reduced-order model. DMD and its improved versions \cite{Chen2012,Williams2015,Rowley2017,Alla2017,Lusch2018}, supported by Koopman operator theory \cite{Budisic2012,Mezic2013}, seek a low-rank approximation to the dynamics of observable data without reliance on the governing equations of motion. With this approach, DMD and Koopman mode expansions are able to linearize the observed dynamics around attracting fixed points on domains that cannot include additional fixed points or limit cycles \cite{Bagheri2013,Brunton2016b,Page2019,Kaiser2021}. Therefore, while truly powerful for globally linearizable dynamics \cite{Page2018}, these linear techniques cannot capture essentially nonlinear dynamical systems (or \textit{non-linearizable} systems) with multiple coexisting steady states.

Other approaches treat the dimensionality reduction and the data-driven dynamical modeling as separate problems. Typically, the data is first processed via a dimensionality reduction algorithm, which ranges from POD or Principal Component Analysis (PCA) \cite{Pearson1901}, its kernelized version \cite{Bishop2006}, manifold learning techniques \cite{Roweis200,Bengio2013,Loiseau2020} or autoencoders \cite{Bengio2016,Champion2019}. Afterwards, the dynamics are identified in the reduced coordinates using classic regression techniques (least-squares \cite{Bishop2006}, LASSO \cite{Tibshirani1996}, SINDy \cite{Brunton2016}), Bayesian learning techniques \cite{Abdessalem2019} or neural networks in different architectures (fully connected, convolutional, recurrent) \cite{Billings2013,Bengio2016,Hartman2017,Lui2019,Raissi2019}. Some of these techniques return complex, black-box models (which may be non-physical \cite{Karniadakis2021}), while others offer sparse models (LASSO, SINDy, Bayesian learning), which allow for easy interpretation and analysis \cite{Brunton2016}. The resulting dynamics, however, are intrinsically determined by the representation offered by dimensionality reduction algorithms, unless penalized in the optimization \cite{Champion2019}. Indeed, the advocated simplicity of sparse models depends critically on the reduction method, as even a linear coordinate change will dramatically destroy the sparsity of a model. In addition, those methods feature a high number of hyperparameters that need to be tuned extensively for good performance. Most importantly, the eventual lack of predictive capabilities often makes the models unattractive for practical use. Indeed, the insertion of parameter variations, disturbances or external forcing into these models is generally heuristic, and hence returns questionable conclusions. 

Our objective here is to discuss a new data-driven reduced-order modeling approach in the context of mechanical vibrations. Based on the recent theory of spectral submanifolds (SSMs) \cite{Haller2016}, this approach identifies very low dimensional, sparse models over different time scales by restricting the full system dynamics to a nested family of attractors. The SSMs forming this family are the smoothest nonlinear continuations of the eigenspaces of the linear part of the dynamical system. When transformed to a normal form, the reduced dynamics on each SSM is low-dimensional, sparse and relevant for all trajectories in the domain of attraction of the SSM. Importantly, each SSM may contain multiple coexisting steady states and hence capture non-linearizable dynamics.

The details and several applications of SSM theory are discussed in \cite{Haller2016,Haller2017,Szalai2017,Jain2018,Ponsioen2018,Breunung2018,Ponsioen2019,Ponsioen2020,Jain2021} and an open-source \textsc{Matlab}\textsuperscript{\textregistered} implementation, \texttt{SSMTool}, for an arbitrary, finite-dimensional dynamical system is available in \cite{SSMTool2021}. Another concept closely linked to SSMs is that of invariant foliations \cite{Szalai2020}, which provides a rigorous nonlinear extension of classic linear modal analysis. Our present discussion of data-driven SSM-based models follows the terminology and notation of the more technical exposition in \cite{Cenedese2021}. The open-source \textsc{Matlab}\textsuperscript{\textregistered} code, \texttt{SSMLearn}, used in our examples is available in \cite{SSMLearn}.

The remainder of this paper is organized as follows. Section \ref{sec:S2} introduces SSMs and discusses their relevance for data-driven model reduction, also depending on the type of experiments that generate the data. We also discuss how our method is complementary to (non-parametric) signal processing techniques in nonlinear system identification \cite{Noel2017}, ranging from the Hilbert transform and its variants \cite{Feldman2011,Jin2020} to wavelet decompositions \cite{Moore2018}. In section \ref{sec:S2}, we summarize our data-driven identification of SSMs and the resulting explicit models on SSMs. We demonstrate the method in section \ref{sec:S3} on a set of examples, which are all analyzed via the \textsc{Matlab}\textsuperscript{\textregistered} implementation of our approach \texttt{SSMLearn} \cite{SSMLearn}. In particular, after a preliminary numerical example, we examine two experimental datasets of nonlinear mechanical systems, one of which regards an internally resonant structure. The data sets in these examples come from diverse sources, from non-contact measurement systems (e.g., digital image correlation, laser scanner vibrometry) to classic accelerometers.

\section{Spectral submanifolds and data-driven models on them}
\label{sec:S2}
In this paper, we consider $N$-degree-of-freedom mechanical systems of the form
\begin{equation}
\label{eq:mechsys}
\mathbf{M}(\mathbf{q})\ddot{\mathbf{q}} = \mathbf{f}(\mathbf{q},\dot{\mathbf{q}}), \,\,\,\,\,\,\,\,\, \mathbf{f}(\mathbf{0},\mathbf{0})=\mathbf{0}, \,\,\,\,\,\,\,\,\, \mathbf{q} \in \mathbb{R}^N, \,\,\,\,\,\,\,\,\, N\geq 1,
\end{equation}
where $\mathbf{q}$ is a generalized coordinate vector, $\mathbf{M}(\mathbf{q})\in\mathbb{R}^{N\times N}$ is a positive definite, symmetric mass matrix. The forcing vector $\mathbf{f}(\mathbf{q},\dot{\mathbf{q}})$ contains all conservative and non-conservative forces, including linear and nonlinear ones. The matrix $\mathbf{M}(\mathbf{q})$, its inverse and $\mathbf{f}(\mathbf{q},\dot{\mathbf{q}})$ are of class $C^r$ with $r\in\mathbb{N}^+ \cup \{\infty\}$ (smooth functions) or $r=a$ (analytic functions).

The equivalent first-order form of Eq. (\ref{eq:mechsys}), with $\mathbf{x} = (\mathbf{q},\dot{\mathbf{q}})\in\mathbb{R}^{n}$ and $n = 2N$, reads
\begin{equation}
\label{eq:dynsys}
\dot{\mathbf{x}} = \mathbf{A}\mathbf{x} + \begin{pmatrix} \mathbf{0} \\ \mathbf{b}(\mathbf{x}) \end{pmatrix}, \,\,\,\,\,\mathbf{A} = \begin{bmatrix} \mathbf{0} & \mathbf{I} \\ \mathbf{M}^{-1}(\mathbf{0})D_{\mathbf{q}}\mathbf{f}(\mathbf{0},\mathbf{0}) & \mathbf{M}^{-1}(\mathbf{0})D_{\dot{\mathbf{q}}}\mathbf{f}(\mathbf{0},\mathbf{0})\end{bmatrix},
\end{equation}
where $\mathbf{b}(\mathbf{x}) = \mathbf{M}^{-1}(\mathbf{q})\mathbf{f}(\mathbf{q},\dot{\mathbf{q}}) - \mathbf{M}^{-1}(\mathbf{0})D_{\mathbf{q}}\mathbf{f}(\mathbf{0},\mathbf{0})\mathbf{q} - \mathbf{M}^{-1}(\mathbf{0})D_{\dot{\mathbf{q}}}\mathbf{f}(\mathbf{0},\mathbf{0})\dot{\mathbf{q}}$. We assume that $\mathbf{x}=\mathbf{0}$ is an asymptotically stable equilibrium and that $\mathbf{A}$ is a semi-simple matrix and has $N$ complex conjugate pairs of eigenvalues with negative real parts. We order these eigenvalues $\lambda_1,\bar{\lambda}_1,\lambda_2,\bar{\lambda}_2,...,\lambda_N,\bar{\lambda}_N$ with decreasing real parts, and we denote by $E_1,E_2,...,E_N$ the corresponding two-dimensional eigenspaces (or modal subspaces).

We denote by $E^{2m}$ the direct sum of $m$ of these modal subspaces, i.e., $E^{2m} = E_{j_1} \oplus E_{j_2} \oplus ... ,\oplus E_{j_m}$. The $2m$-dimensional, spectral subspace $E^{2m}$ is invariant for the linearization of system (\ref{eq:dynsys}). Its reduced dynamics is governed by the eigenvalues $\lambda_{j_1},\lambda_{j_2},...,\lambda_{j_m}$, which, along with the conjugate ones, form the set $\mathrm{Spect}\left(\mathbf{A}|_{E^{2m}}\right)$. The spectral submanifold (SSM), $\mathcal{W}(E^{2m})$, is the smoothest nonlinear continuation of the linear subspace $E^{2m}$ \cite{Haller2016}, as can be deduced from the more abstract invariant manifold results of \cite{Cabre2003a,Cabre2003b,Cabre2005,Haro2016}. Specifically, $\mathcal{W}(E^{2m})$ is the unique $2m$-dimensional, class $C^r$ invariant manifold of system (\ref{eq:dynsys}) tangent to the spectral subspace $E^{2m}$ at the origin. The existence of $\mathcal{W}(E^{2m})$ is guaranteed whenever the eigenvalues $(\lambda_j,\bar{\lambda}_j)$ with $j \neq j_1,j_2,...,j_m$ are not in resonance with those in $\mathrm{Spect}\left(\mathbf{A}|_{E^{2m}}\right)$, i.e., for $\mathbf{k} = (k_1,k_2,...,k_{2m})\in\mathbb{N}^{2m}$,
\begin{equation}\label{eq:resonance}
\lambda_j - \sum_{l=1}^{m} \big( \lambda_{j_l} k_l  + \bar{\lambda}_{j_l}k_{l+m} \big) \neq 0, \,\,\,\,\,\,\,\, \sum_{l=1}^{2m}k_l \leq \mathrm{Int}\left[\frac{\displaystyle \min_{\lambda\in\mathrm{Spect}(\mathbf{A})-\mathrm{Spect}\left(\mathbf{A}|_{E^{2m}}\right)}\mathrm{Re}(\lambda)}{\displaystyle \max_{\lambda\in\mathrm{Spect}\left(\mathbf{A}|_{E^{2m}}\right)}\mathrm{Re}(\lambda)}\right],
\end{equation}
as discussed in \cite{Haller2016}. From a numerical perspective, the nonresonance condition in Eq. (\ref{eq:resonance}) is violated if the absolute value of the left-hand side of the inequality is below a certain tolerance. In that case, one needs to add the resonant modal subspace $E_j$ to $E^{2m}$, resulting in the SSM of the form $\mathcal{W}(E^{2m}\oplus E_j)$. This larger SSM can be used to capture nonlinear modal interactions, e.g., in weakly damped systems with rationally dependent linearized frequencies. The dynamics restricted to SSMs gives exact nonlinear reduced-order models for system (\ref{eq:dynsys}) \cite{Haller2016,Ponsioen2018,Jain2021}.
\begin{figure}[t]
    \centering
    \includegraphics[width=1\textwidth]{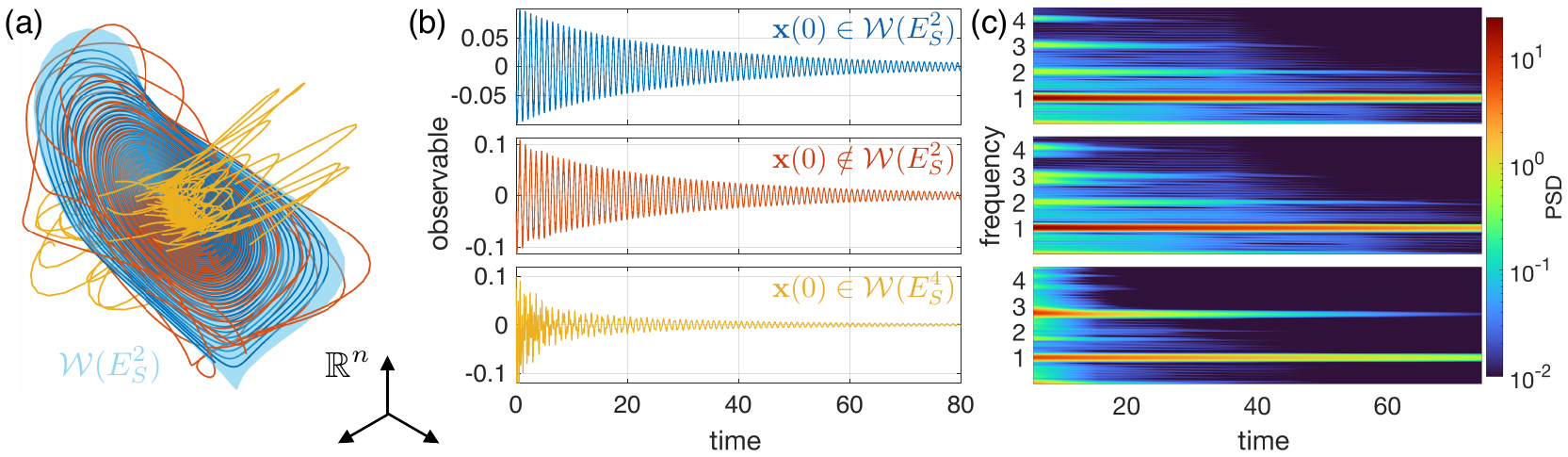}
    \caption{Illustration of different time scales of nonlinear dynamics captured by a nested set of SSMs. (a) Dynamics near a slow two-dimensional SSM $\mathcal{W}(E^{2}_S)$ with three trajectories in the phase space. (b,c) The three trajectories, shown with consistent colors, and their spectrogram. These trajectories were generated by a mechanical system in which $1, 2.7, 4.2$ are the first three linearized frequencies. In (b,c) the top trajectory has initial condition $\mathbf{x}(0)$ on the slow two-dimensional SSM $\mathcal{W}(E^{2}_S)$ and decays on it, the middle is initialized with a small perturbation off $\mathcal{W}(E^{2}_S)$, and the bottom one decays on $\mathcal{W}(E^{4}_S)$.}
    \label{fig:S2_IM1}
\end{figure}

The most important SSMs from a data-driven perspective are slow SSMs, which are constructed over the spectral subspace spanned by the $m$ slowest modes: $E^{2m}_S = E_{1} \oplus E_{2} \oplus ... \oplus E_{m}$. Slow SSMs are attracting normally hyperbolic invariant manifolds \cite{Haller2016} to which nearby trajectories converge exponentially fast \cite{Fenichel1974}, as illustrated in Fig. \ref{fig:S2_IM1}(a). Therefore, generic experiments on mechanical systems in the form of (\ref{eq:dynsys}) will yield trajectories converging exponentially fast to slow SSMs, which in turn capture the asymptotic dynamics near the equilibrium. Faster time-scales of the dynamics can be extracted from trajectory data by model reduction to higher dimensional members of the nested slow SSM family $\mathcal{W}(E^{2}_S)\subset\mathcal{W}(E^{4}_S) \subset ... \subset \mathcal{W}(E^{2(n-1)}_S) \subset \mathbb{R}^n$. As an illustration, trajectories on slow SSMs are shown in Fig. \ref{fig:S2_IM1}(b), along with their spectrogram (or short time Fourier transform) in Fig. \ref{fig:S2_IM1}(c). The trajectory on the top is in $\mathcal{W}(E^{2}_S)$, while the middle one is initialized close to $\mathcal{W}(E^{2}_S)$ and hence converges to $\mathcal{W}(E^{2}_S)$, as seen from the disappearance of higher frequencies. Finally, the trajectory at the bottom in Fig. \ref{fig:S2_IM1}(b,c) belongs to $\mathcal{W}(E^{4}_S)$, where two modal contributions can be clearly identified.

For trajectories with generic initial conditions---as those generated by hammer impacts---discarding the initial part of the measured signal yields trajectory data close to a slow $2m$-dimensional SSM, where $m$ is the number of dominant frequencies in the signal. This number $m$ is inferred from a preliminary time-frequency analysis (e.g. spectrograms, wavelet transforms \cite{Moore2018}), as those of Fig. \ref{fig:S2_IM1}(c). The larger $m$, the more data is needed to properly explore the SSM. In contrast, targeted experiments can focus on specific SSMs, which is the case for two-dimensional SSMs in resonance decay experiments \cite{Peeters2011a,Ehrhardt2016}. In that setting, a near-resonant oscillation is first isolated using a shaker, which is then turned off. This generates a system trajectory that decays towards the equilibrium along the targeted two-dimensional SSM, provided that no internal resonance occurs. Using the shaker, we can typically excite higher amplitudes in comparison to hammer impact testing, where the energy is spread among multiple modes.

\subsection{Learning SSMs from data}
\label{sec:S2_1}
To learn SSMs from data, we use the methodology presented in \cite{Cenedese2021}, which is implemented in the open-source \textsc{Matlab}\textsuperscript{\textregistered} package, \texttt{SSMLearn} \cite{SSMLearn}. In what follows, we sketch the main ideas of this method before going into the details of the data-driven reduced-order models that \texttt{SSMLearn} can identify.

Measuring all phase space variables of a mechanical system is generally unrealistic. Typically, only a limited set of observed quantities is available, so that we need to embed the SSM, $\mathcal{W}(E^{2m})$, into a lower dimensional space of observables. According to the prevalence version of Whitney's embedding theorem \cite{Sauer1991}, almost all sets of independent and simultaneous measurements $\mathbf{y}(t) = (y_1(t),y_2(t),...y_p(t))\in\mathbb{R}^p$ form an embedding space for $2m$-dimensional SSMs if $p>4m$. This is the case, for example, when displacements and velocities of at least $3m$ material points of a mechanical system are available. Practical experiments, however, generally only provide the displacement, velocity or acceleration of a single material point, denoted as $s(t)\in\mathbb{R}$, recorded at $\Delta t$ time intervals. To this end, we exploit Takens delay embedding theorem \cite{Takens1981}, which, in its prevalence version \cite{Sauer1991}, states that $\mathbf{y}(t) = (s(t),s(t+\Delta t),s(t+2\Delta t),...s(t+(p-1)\Delta t))\in\mathbb{R}^p$ forms an embedding space with probability one if $p>4m$ under generic nondegeneracy conditions on the sampling time $\Delta t$. Further spaces may also qualify in practice, e.g., featuring $p\leq 4m$ or constructed from multiple measurements augmented by delays, but one needs to examine on a case-by-case basis whether these are embedding spaces or not.

We denote by $\mathcal{M}_0$ the embedded SSM, for which now we need to construct a reduced-order model in the embedding space. We assume that the equilibrium is at $\mathbf{y}=\mathbf{0}$ and that $\mathcal{M}_0$ does not fold over its tangent space at the origin $T_{\mathbf{0}}\mathcal{M}_0$, so that we can construct a data-driven graph-style parametrization for $\mathcal{M}_0$ over $T_{\mathbf{0}}\mathcal{M}_0$. We let $\mathbf{V}_1\in\mathbb{R}^{p\times 2m}$ be the matrix whose orthonormal columns span $T_{\mathbf{0}}\mathcal{M}_0$ and we define the SSM parametrization, $\mathbf{v}:\mathbb{R}^{2m}\rightarrow\mathbb{R}^{p}$, as
\begin{equation}\label{eq:SSMpara}
\mathbf{y}=\mathbf{v}(\mathbf{V}_1^\top\mathbf{y})=\mathbf{V}_1\mathbf{V}_1^\top\mathbf{y} + \mathbf{v}_{\mathrm{nl}}(\mathbf{V}_1^\top\mathbf{y}), \,\,\,\,\,\,\,\, \mathbf{V}_1^\top \mathbf{V}_1 = \mathbf{I}, \,\,\,\,\,\,\,\, \mathbf{V}_1^\top \mathbf{v}_{\mathrm{nl}}(\mathbf{V}_1^\top\mathbf{y}) = \mathbf{0},
\end{equation} 
where we assume that $\mathbf{v}_{\mathrm{nl}}:\mathbb{R}^{2k}\rightarrow \mathbb{R}^{p}$ is a multivariate polynomial from order $2$ to $M$. The matrix $\mathbf{V}_1$, as well as the coefficients of the polynomial $\mathbf{v}_{\mathrm{nl}}$, can be found via constrained maximum likelihood estimation of (\ref{eq:SSMpara}), as discussed in \cite{Cenedese2021}.

Once trajectories in the projection coordinates $\mathbf{V}_1^\top\mathbf{y}\in\mathbb{R}^{2m}$ are known, we can identify the SSM-reduced dynamics. Here, the idea is to find the extended normal form of the vector field governing the dynamics in the projection coordinate (or reduced) domain \cite{Cenedese2021}, motivated by classic studies of bifurcations \cite{GH1983,Murdock2003}. Specifically, we need to find an invertible change of coordinates $\mathbf{V}_1^\top\mathbf{y}=\mathbf{h}(\mathbf{z})$ (and its inverse) that brings the SSM-reduced dynamics to its simplest possible complex polynomial form $\dot{\mathbf{z}}=\mathbf{n}(\mathbf{z})$ with $\mathbf{z}\in\mathbb{C}^{2m}$. The linear part of $\mathbf{n}$ is the diagonal matrix of the eigenvalues related to the SSM, with $\mathbf{z}=(z_1,\bar{z}_1,z_2,\bar{z}_2,...z_m,\bar{z}_m)$ denoting complex modal coordinates for the linearized system. The maps $\mathbf{h}$, $\mathbf{h}^{-1}$ and $\mathbf{n}$ are multivariate polynomials with their coefficients determined from an extended normal form approach used in classic unfoldings of bifurcations \cite{GH1983,Murdock2003}. In this approach, the classic Poincaré \cite{Poincare1892} normal form construct is relaxed in that not only resonant but also near-resonant terms are kept in the normal form (see \cite{Cenedese2021} for more details). This normalization renders $\mathbf{n}$ a sparse vector field extracting the fundamental physics, as we discuss in the next section. We determine resonant coefficients from an initial estimate of the linearized dynamics, and we identify from data the maps $\mathbf{h}$, $\mathbf{h}^{-1}$ and $\mathbf{n}$ by minimizing the conjugacy error, as explained in detail in \cite{Cenedese2021}. For example, the structure of the cubic normal form for a two-dimensional SSM is
\begin{equation}\label{eq:2SSMcdyn}
\begin{array}{c}
\mathbf{z} = (z,\bar{z}), \,\,\,\,\,\,\,\, \mathbf{h}(\mathbf{z}) = (h_1,\bar{h}_1), \,\,\,\,\,\,\,\, \mathbf{n}(\mathbf{z}) = (n_1,\bar{n}_1), \\
h_1(\mathbf{z}) = z + h_{20}z^2 + h_{11}z\bar{z} + h_{02}\bar{z}^2  + h_{30}z^3 + h_{12}z\bar{z}^2 + h_{03}\bar{z}^3, \,\,\,\,\,\,\,\,
n_1(\mathbf{z}) = \lambda z + \gamma z^2\bar{z},
\end{array}
\end{equation} 
which resembles the classic Hopf normal form \cite{Marsden1976}. These normal form models are particularly simple to handle in polar coordinates $(\rho_j,\theta_j)$, defined as $z_j = \rho_j e^{i\theta_j}$ for $j = 1,2,..., m$.

\subsection{Interpretability and extrapolation from SSM-reduced models}
\label{sec:S2_2}
The most general normal form on a $2m$-dimensional SSM is
\begin{equation}\label{eq:SSMdyn}
\begin{array}{l}
\dot{\rho}_j = -\alpha_j(\boldsymbol{\rho},\boldsymbol{\theta})\rho_j, \\
\dot{\theta}_j = \omega_j(\boldsymbol{\rho},\boldsymbol{\theta}),
\end{array}\,\,\,\,\,\,\,\, j = 1,2,..., m, \,\,\,\,\,\,\,\, \boldsymbol{\rho} = (\rho_1,\rho_2,...\rho_m), \,\,\,\,\,\,\,\, \boldsymbol{\theta} = (\theta_1,\theta_2,...\theta_m).
\end{equation} 
Some explicit examples are presented in the examples of section \ref{sec:S3}, including cubic polar normal forms of two-dimensional and four dimensional SSMs, the latter appearing both for non-resonant eigenvalues and for a $1:2$ resonance. If the linearized frequencies are non-resonant, then $\alpha_j$ and $\omega_j$ only depend on the amplitudes $\boldsymbol{\rho}$. The normal form (\ref{eq:SSMdyn}) then decouples the amplitude dynamics from the phase dynamics. This enables us to distinguish different modal contributions, perform a slow-fast decomposition, detect modal interactions and analyze the uncoupled oscillator limit. The zero-amplitude limit of the functions $\alpha_j$ and $\omega_j$ converges to the linearized damping and frequency of mode $j$, i.e.,
\begin{equation}
\lim_{\|\boldsymbol{\rho}\|\rightarrow 0} \big[ -\alpha_j(\boldsymbol{\rho},\boldsymbol{\theta}) + i \omega_j(\boldsymbol{\rho},\boldsymbol{\theta}) \big]= \lambda_j.
\end{equation} 
Hence, $\alpha_j$ and $\omega_j$ are the nonlinear continuations of these linear quantities, characterizing how dissipation and frequency change with respect to the amplitudes (and phases for internally resonant systems). For a two-dimensional SSM, the parametrized curves $\alpha(\rho)$ and $\omega(\rho)$ are the backbones of transient oscillations \cite{Peeters2011a,Szalai2017,Breunung2018,TRC2019_PartII}, representing the instantaneous damping and frequency as nonlinear functions of the normal form amplitude $\rho$. Normal form amplitudes do not, however, have any direct physical meaning. For physical insights, we need to express any amplitude of interest via the SSM parametrization $\mathbf{v}$ and the normal form transformation $\mathbf{h}$. For instance, for two-dimensional SSMs and for a scalar quantity $g:\mathbb{R}^p\rightarrow \mathbb{R}$ defined on the observable space $\mathbb{R}^p$, the amplitude of the oscillations can be defined as \cite{Szalai2017,Ponsioen2018}
\begin{equation}\label{eq:NFphysamp}
\mathrm{amp}(\rho) = \max_{\theta \in [0,2\pi)} \left| g \left( \mathbf{v} \left( \mathbf{h} \left( \mathbf{z} \right) \right)\right)\right|, \,\,\,\,\,\,\,\, \mathbf{z} = \left( \rho e^{i\theta},\rho e^{-i\theta} \right).
\end{equation} 
Then, backbone curves can be expressed as parametric curves $\{\alpha(\rho),\mathrm{amp}(\rho)\}$ and $\{\omega(\rho),\mathrm{amp}(\rho)\}$. 

SSMs are robust features of the dynamics, because they survive under small autonomous perturbations and even under some non-autonomous perturbations of the vector field (\ref{eq:dynsys}) \cite{Haller2016}. The most important class of these perturbations in our context is that of small external time-periodic forcing appearing on the right-hand side of Eq. (\ref{eq:mechsys}). In that case, the autonomous SSM will serve as the leading order approximation for a non-autonomous, time-periodic SSM that carries reduced, time-periodic dynamics \cite{Haller2016,Breunung2018,Ponsioen2019,Jain2021}. With the addition of such forcing, the normal form (\ref{eq:SSMdyn}) becomes \cite{Cenedese2021}
\begin{equation}\label{eq:SSMdynf}
\dot{\rho}_j = -\alpha_j(\boldsymbol{\rho},\boldsymbol{\theta})\rho_j - f_j\sin(\Omega t - \theta_j), \,\,\,\,\,\,\,\, \dot{\theta}_j = \omega_j(\boldsymbol{\rho},\boldsymbol{\theta}) + \frac{f_j}{\rho_j}\cos(\Omega t - \theta_j),
\end{equation} 
where $\Omega$ is the forcing frequency and $f_j$ the forcing amplitudes for each mode. Generally, numerical continuation is necessary for studying periodic responses and eventual bifurcations of (\ref{eq:SSMdynf}) depending on forcing frequencies and amplitudes. In the simplest case of $m=1$, however, we can introduce the phase shift $\psi = \theta-\Omega t$ to obtain from (\ref{eq:SSMdynf}) the forced normal form
\begin{equation}\label{eq:2SSMdynf}
\dot{\rho} = -\alpha(\rho)\rho + f\sin(\psi), \,\,\,\,\,\,\,\, \dot{\psi} = \omega(\rho)-\Omega + \frac{f}{\rho}\cos(\psi),
\end{equation} 
which yields closed form predictions for amplitudes and phases of the forced periodic solutions:
\begin{equation}\label{eq:2SSMFRC}
\Omega = \omega(\rho)\pm\sqrt{\frac{f^2}{\rho^2}-\alpha^2(\rho)}, \,\,\,\,\,\,\,\,\psi = \sin^{-1}\left(\frac{\alpha(\rho)\rho}{f} \right),
\end{equation} 
known as frequency response curves (FRCs), parametrized by the amplitude $\rho$. Predictions of these curves from unforced data, however, have generally been unavailable. Physical amplitudes can be derived from the predictions of (\ref{eq:2SSMFRC}) using Eq. (\ref{eq:NFphysamp}) and the stability of the predicted forced response can be derived from the Jacobian of the vector field (\ref{eq:2SSMdynf}) \cite{Breunung2018}. We note find from Eq. (\ref{eq:2SSMFRC}) that the forced backbone curve (the location of maximal amplitude responses of FRCs under varying $f$) coincides with that of decaying oscillations, given by $\omega(\rho)$. Specifically, maximal amplitude responses occur at amplitudes $\rho_{\max}$ satisfying $f = \alpha(\rho_{\max})\rho_{\max}$, $\Omega = \omega(\rho_{\max})$ and phase-lag quadrature, i.e., $\theta = \Omega t -\pi/2$. These maximal amplitude responses can be used to calibrate the normal form forcing amplitude $f$ with experimentally exerted forcing levels.

Equations (\ref{eq:SSMdynf},\ref{eq:2SSMdynf}) have $O(f\rho)$ accuracy \cite{Breunung2018,Cenedese2021}, but higher-order approximations can improve this accuracy further \cite{Ponsioen2020}. We expect, for example, that forced backbone curves depart from those of decaying oscillations at large motion and/or large forcing amplitude values \cite{Cenedese2019}. From a data-driven perspective, once the autonomous core of Eqs. (\ref{eq:SSMdynf},\ref{eq:2SSMdynf}) is identified, we only need to calibrate the forcing amplitudes for predicting FRCs. In matching experimental results, one calibration point is sufficient if the forcing amplitude is kept constant during experimental frequency sweeps. The SSM-parametrization under forcing is also modulated by a small time-periodic component \cite{Ponsioen2020,Cenedese2021}. This may be disregarded since the leading-order non-linearizable dynamics is already captured by forcing in the reduced dynamics (\ref{eq:SSMdynf},\ref{eq:2SSMdynf}).
\begin{figure}[t!]
    \centering
    \includegraphics[width=1\textwidth]{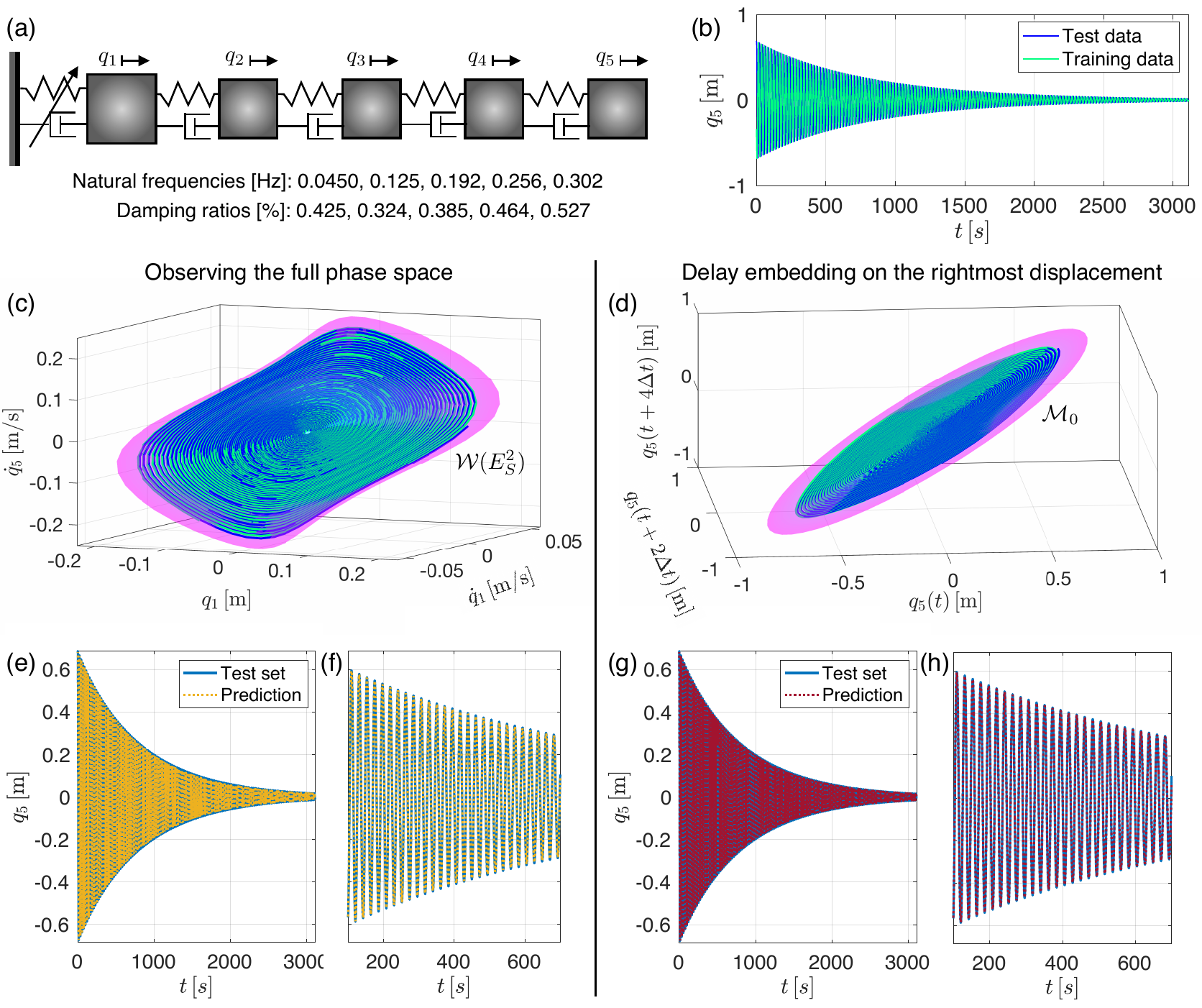}
    \caption{(a) Sketch of the oscillator chain considered in Sec. \ref{sec:S3}\ref{sec:S3_OC}. (b) Two trajectories decaying on its slow two-dimensional SSM. (c,d) The SSM and the trajectories in the phase space and in the delay observable space, respectively. (e-h) The performance of the normalized SSM-reduced models in reconstructing the test trajectory.}
    \label{fig:S3_OC1}
\end{figure}

\section{Examples}
\label{sec:S3}
We now discuss some examples that illustrate the power of the SSM-based, data-driven model reduction method we have discussed. Our first example is a chain of lumped oscillators, while the other two involve data from laboratory experiments. Additional details and further examples can also be found in the \textsc{Matlab}\textsuperscript{\textregistered} live-scripts in \cite{SSMLearn}.

To express trajectory reconstruction errors, we use the normalized mean trajectory error NMTE, which, for a dataset of $P$ instances of observable points $\mathbf{y}_j\in\mathbb{R}^p$ and their reconstruction $\hat{\mathbf{y}}$, is defined as
\begin{equation}
\mathrm{NMTE} = \frac{1}{P \| \underline{\mathbf{y}} \|}\sum_{j=1}^P \| \mathbf{y}_j - \hat{\mathbf{y}}_j\|.
\end{equation} 
Here $\underline{\mathbf{y}}$ is a relevant normalization vector, which is usually taken to be the data point $\mathbf{y}_j$ with the maximum norm in the dataset. To validate the reduced dynamics on a test trajectory, we integrate the reduced-order model from the same initial condition and compare the results. 

\subsection{Identification of SSMs in a chain of oscillators}
\label{sec:S3_OC}
We consider the chain of oscillators sketched in Fig. \ref{fig:S3_OC1}(a), where we set the first mass as $1.5$ [kg] and the others as $1$ [kg]. We also assume all spring-dampers to be linear with unitary stiffness, except for the leftmost one that exerts a nonlinear force $f_{nl,1} = 0.33\dot{q}_1^2 + 3 q_1^3 + 0.7q_1^2\dot{q} + 0.5\dot{q}_1^3$ on the first mass. The linear damping matrix for the system is proportional to the mass and stiffness matrices with constants $0.002$ and $0.005$, with the resulting eigenvalues at the trivial equilibrium reported in Fig. \ref{fig:S3_OC1}(a).

\begin{figure}[t!]
    \centering
    \includegraphics[width=1\textwidth]{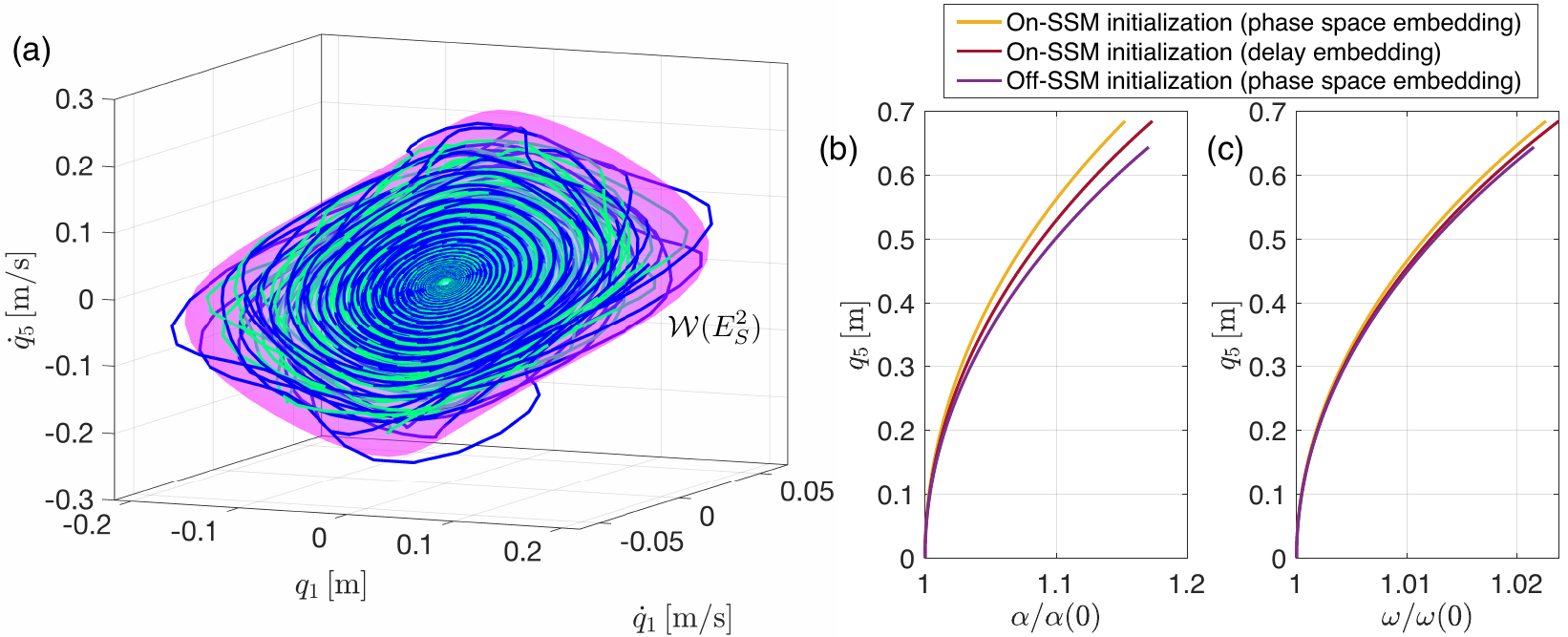}
    \caption{(a) Two trajectories converging to the slowest two-dimensional SSM of the oscillator chain. (b,c) Instantaneous damping and frequency curves constructed from a phase space embedding with perfect SSM initialization, from a delay embedding with perfect SSM initialization and from a phase space embedding with imperfect SSM initialization.}
    \label{fig:S3_OC2}
\end{figure}

\begin{figure}[t!]
    \centering
    \includegraphics[width=1\textwidth]{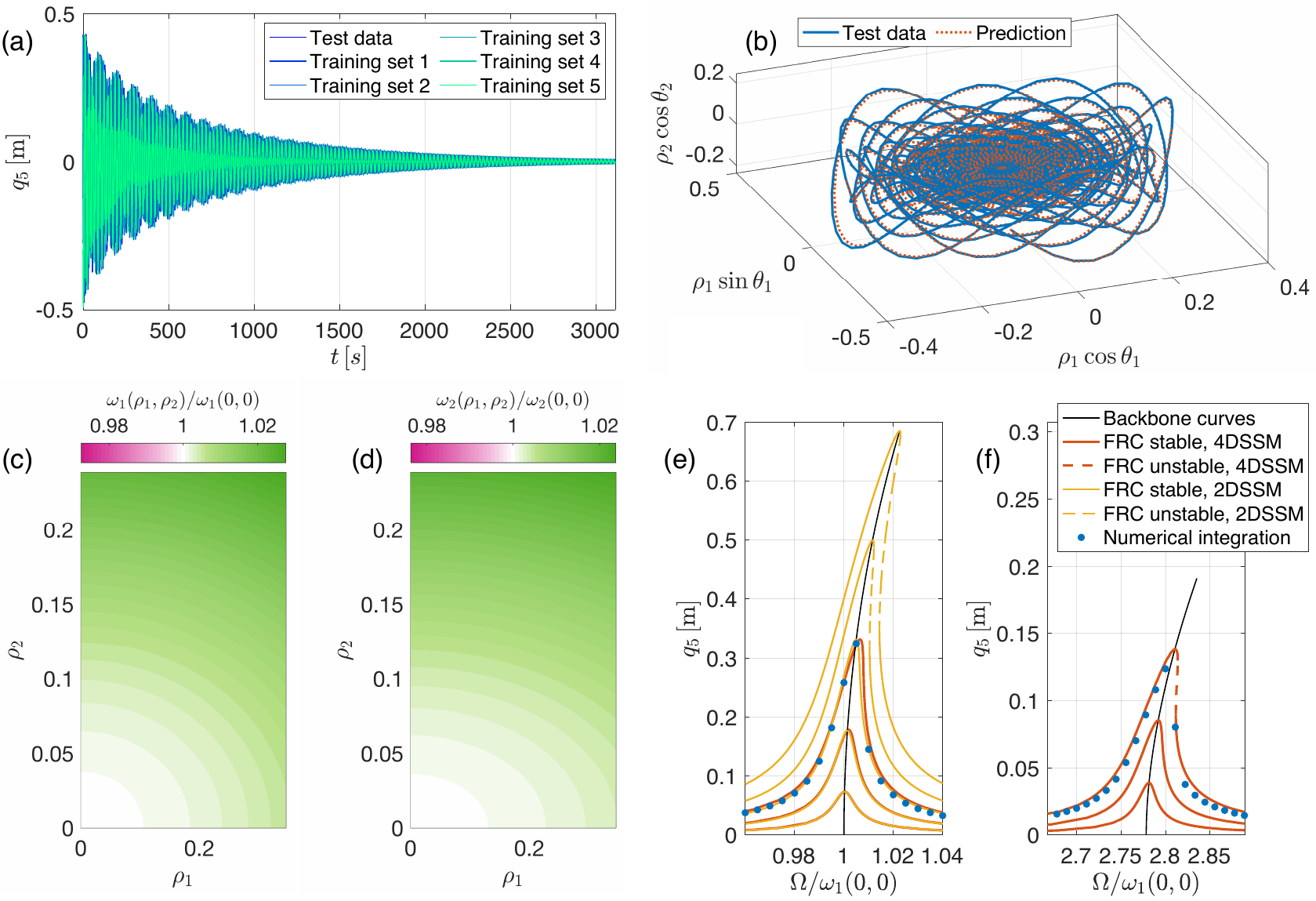}
    \caption{(a) Decaying trajectories from the slow four-dimensional SSM $\mathcal{W}(E_S^4)$ of the oscillator chain. (b) Test trajectory and its model-based prediction in the normal form domain. (c,d) Instantaneous frequencies of the slow (mode $1$) and fast ($2$) modes of $\mathcal{W}(E_S^4)$. (e,f) Frequency Response Curve (FRC) and backbone curves predictions from the reduced-order models (\ref{eq:RDOC2}) and (\ref{eq:RDOC4}) along with forced steady states (dots) obtained via numerical integration of the full system.}
    \label{fig:S3_OC3}
\end{figure}

We start with the study of the slow two-dimensional SSM $\mathcal{W}(E^{2}_S)$ of the oscillator chain. We compute this SSM via \texttt{SSMTool} \cite{SSMTool2021}, from which we initialize the two decaying trajectories shown in Fig. \ref{fig:S3_OC1}(b); one of these trajectories is used for testing the constructed model. We identify reduced-order models from two different observables. The first observable set is the set of all phase space variables, while the second is a set of delayed samples of the (scalar) displacement of the rightmost mass $q_5$. We select the delay embedding of minimal dimension (five) required by the Takens theorem. The embedding coordinates are, therefore, $\mathbf{y}(t) = (q_5(t),q_5(t+\Delta t),q_5(t+2\Delta t),q_5(t+3\Delta t),q_5(t+4\Delta t))$, where the sampling time $\Delta t$ is $0.445$ [s]. A cubic-order parametrization for the phase space embedding and a parametrization for the delay embedding show good accuracy. The SSM $\mathcal{W}(E^{2}_S)$ and its embedding in the delay space $\mathcal{M}_0$ are shown in Fig. \ref{fig:S3_OC1}(b,c). The flat appearance of the manifold in \ref{fig:S3_OC1}(d) in the delay space is a general phenomenon, as shown mathematically in \cite{Cenedese2021}. The cubic polar normal form on the phase-space-embedded SSM is found by \texttt{SSMLearn} to be
\begin{equation}\label{eq:RDOC2}
\dot{\rho} = -0.001201\rho-0.0007300\rho^{3} = -\alpha(\rho)\rho , \,\,\,\,\,\,\,\, 
\dot{\theta} = +0.2827+0.02546 \rho^{2} = \omega(\rho).
\end{equation} 
A similar model is identified for the delay embedding. Both reduced-order models capture well the dynamics of the testing trajectories, as seen in Fig. \ref{fig:S3_OC1}(e-h), with less than $2$ $\%$ NMTE error. The instantaneous damping $\alpha(\rho)$ and frequency $\omega(\rho)$ are shown in \ref{fig:S3_OC2}(b,c), displaying only a minimal disagreement. We note that this identification is robust against perturbations of the initial condition. Indeed, if we initialize trajectories slightly off the SSM as shown in Fig. \ref{fig:S3_OC2}(a), then \texttt{SSMLearn} still finds a good approximation for the reduced dynamics, as demonstrated by the curves in Fig. \ref{fig:S3_OC2}(c). 
If these perturbations are not small enough for the dynamics to be described by a two-dimensional SSMs, then we need to increase the SSM dimension. For instance, we computed trajectories, shown in Fig. \ref{fig:S3_OC3}(a), decaying along the slow four-dimensional SSM, $\mathcal{W}(E^{4}_S)$. Five of these trajectories are used for training and one is left for testing our reduced-order model. 
The cubic normalized, SSM-reduced dynamics identified by \texttt{SSMLearn} has a $2.65$ $\%$ NMTE error and is of the form
\begin{equation}\label{eq:RDOC4}\arraycolsep=1.4pt
\begin{array}{rl}
\dot{\rho}_1 &= -0.001200\rho_1-0.0005548\rho_1^{3}-0.01010\rho_1\rho_2^{2} = -\alpha_1(\rho_1,\rho_2)\rho_1 ,\\ 
\dot{\rho}_2 &= -0.002541\rho_2+0.003728\rho_1\rho_2^{2}-0.05627\rho_2^{3}  = -\alpha_2(\rho_1,\rho_2)\rho_2 ,\\
\dot{\theta}_1 &= +0.2825+0.01316\rho_1^{2}+0.1085\rho_2^{2}  = \omega_1(\rho_1,\rho_2), \\
\dot{\theta}_2 &= +0.7850+0.02340\rho_1^{2}+0.2760\rho_2^{2} = \omega_2(\rho_1,\rho_2). \end{array}
\end{equation}
Eventual differences between the dynamics of the slowest mode in (\ref{eq:RDOC2}) with respect to those in (\ref{eq:RDOC4}) are due to different amplitude scalings. Prediction of a test trajectory based on the model (\ref{eq:RDOC4}) is shown in Fig. \ref{fig:S3_OC3}(b). The instantaneous frequencies for the two modes are shown in Fig. \ref{fig:S3_OC3}(c,d). These are surfaces since both frequencies depend (either weakly or strongly) on both modal amplitudes.

For additional validation, we show FRCs of the models (\ref{eq:RDOC2})-(\ref{eq:RDOC4}) for different forcing amplitudes in Fig. \ref{fig:S3_OC3}(e,f) around the first two eigenfrequencies. While on the two-dimensional SSM $\mathcal{W}(E_S^2)$ we have the closed-form solution (\ref{eq:2SSMFRC}), FRCs on the four-dimensional SSM $\mathcal{W}(E_S^4)$ are computed using the periodic orbit toolbox of the numerical continuation core \textsc{coco} \cite{Dankowicz2013}. These plots in Fig. \ref{fig:S3_OC3}(e,f) are completed with backbones curves and forced responses obtained via numerical integration of the full model. The forcing only acts in the direction of the first two modes, with amplitudes $0.38$ and $1.75$ [mN]. Our data-driven predictions, which are based only on unforced data and a simple calibration procedure for the normal form forcing amplitudes, are in close agreement with the responses from the full system.

\subsection{Resonance decay in the Brake-Reuss beam}
\label{sec:S3_BRB}
\begin{figure}[t!]
    \centering
    \includegraphics[width=1\textwidth]{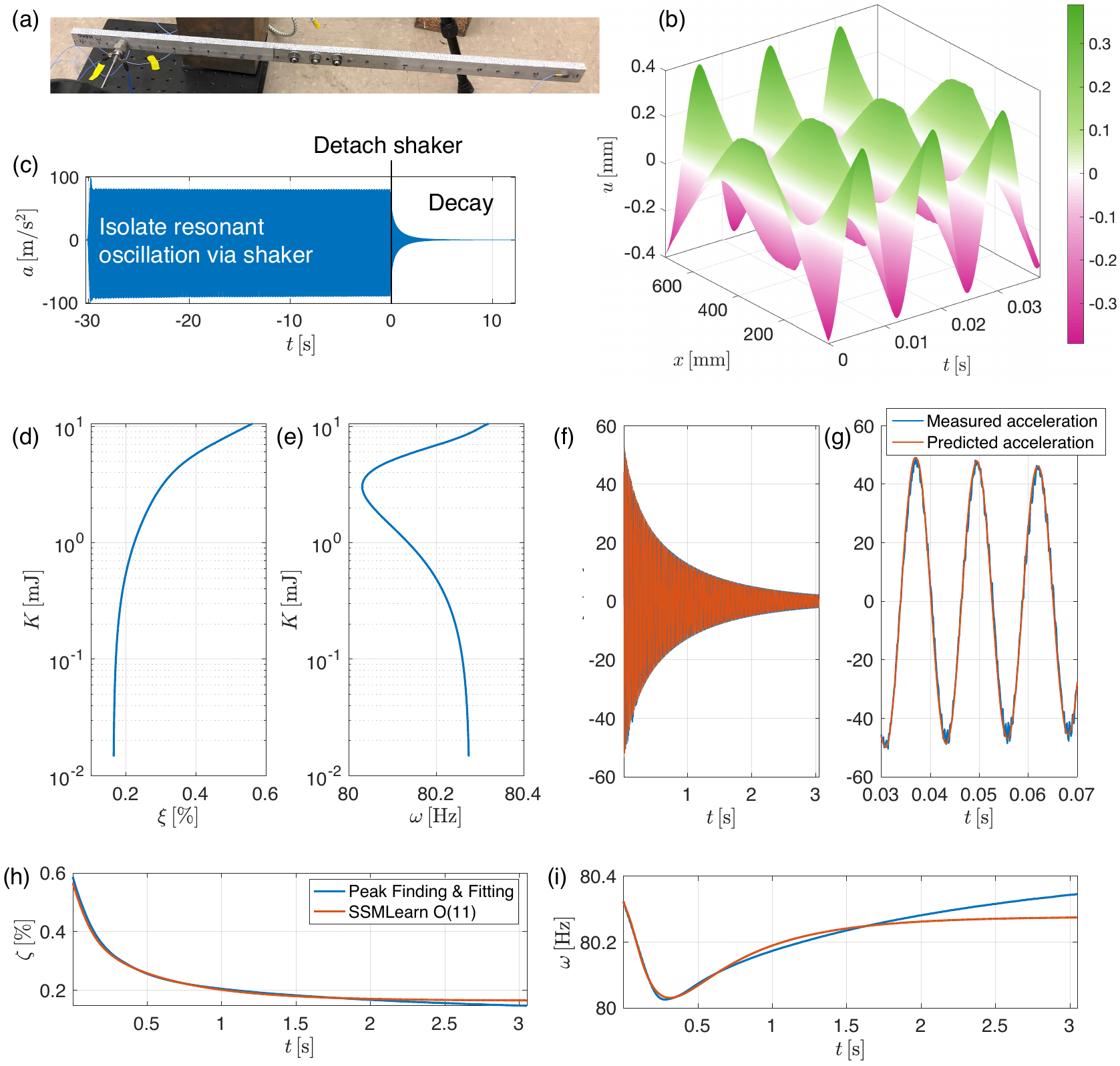}
    \caption{(a) Testing setup for the Brake-Reuss beam. (b,c) The measured displacement and acceleration data. (d-i) Results from the reduced-order model trained on displacement data. The backbone curves in (d,e) show the instantaneous characteristics of the beam with respect to its kinetic energy, while (f,g) validate the predictions of acceleration. Plots (h,i) compare the instantaneous properties of the data-driven model with respect to those measured with the Peak Fitting and Finding method \cite{Jin2020} on the acceleration signal.}
    \label{fig:S3_BRB}
\end{figure}

The Brake-Reuss beam (BRB) is a benchmark system in the study of jointed structures \cite{Brake2018,TRC2019_PartI,TRC2019_PartII}. In our study, it consists of two 304 stainless steel beams assembled with a lap joint, as shown in Fig. \ref{fig:S3_BRB}(a). While full models for these structures may not be smooth, we find that trajectory data can be fitted well to smooth models, thereby justifying an SSM-based approach. The data considered here arises from a single resonance decay test, available from \cite{TRC2019_PartII}, targeting the slowest structural mode. One observable is the measurement from an accelerometer mounted, as shown in Fig. \ref{fig:S3_BRB}(a), on the shaker with time history illustrated in Fig. \ref{fig:S3_BRB}(c). Another observable is the displacement field of the bottom side of the beam, measured using digital image correlation (DIC). The latter dataset, consisting of $206$ points over $72$ [cm] of beam length, has a limited time range due to limitations in camera memory. The initial evolution of the measured displacement field is depicted in Fig. \ref{fig:S3_BRB}(b). 

Our goal in this example is to construct a nonlinear reduced-order model using displacement data and validate it on the acceleration measurement. We truncate the time signals after shaker release to eliminate the influence of disturbances from non-perfect detachment. Nevertheless, high frequency contributions decay rapidly and the transient settles along the slowest SSM. To diversify the data, we augment the displacement with four delayed measurements, so that the observable phase space has dimension $1030$. The SSM is approximately a plane in this space, but the reduced dynamics is highly nonlinear. For adequate accuracy, the normal form indeed needs terms up to $O(11)$ to capture the dynamics:
\begin{equation}\label{eq:RDBRB}
\begin{array}{l}\dot{\rho} = -0.8255\rho-16.05 \rho^{3}+166.3 \rho^{5}-1421 \rho^{7}+5314 \rho^{9}-7138 \rho^{11}=-\alpha(\rho)\rho,\\ \dot{\theta} =+504.4-46.16 \rho^{2}+350.3 \rho^{4}+412.9 \rho^{6}-8468 \rho^{8}+16975 \rho^{10}=\omega(\rho).
\end{array}
\end{equation}
The model can be used to approximate the beam kinetic energy as
\begin{equation}
K = \frac{1}{2}\frac{ m_{\mathrm{BRB}} }{ N_{\mathrm{DIC}} } \sum_{j = 1}^{ N_{\mathrm{DIC}} } \dot{u}_j^2(t),
\end{equation}
where $N_{\mathrm{DIC}} = 206$ is the number of DIC measurement locations and $m_{\mathrm{BRB}} = 1.796$ [kg] is the beam mass. As discussed in \cite{TRC2019_PartII}, the kinetic energy amplitude is a good proxy for the instantaneous decay properties, i.e., the instantaneous damping ratio and frequency, shown in Fig. \ref{fig:S3_BRB}(d,e), respectively. The instantaneous damping ratio is defined from the normal form dynamics (\ref{eq:RDBRB}) as $\xi(\rho) = \alpha(\rho)/\omega(\rho)$, expressed in percentage. The damping exhibits a strong variation from its linear limit, while the frequency here shows a peculiar softening-hardening trend.
\begin{figure}[t!]
    \centering
    \includegraphics[width=1\textwidth]{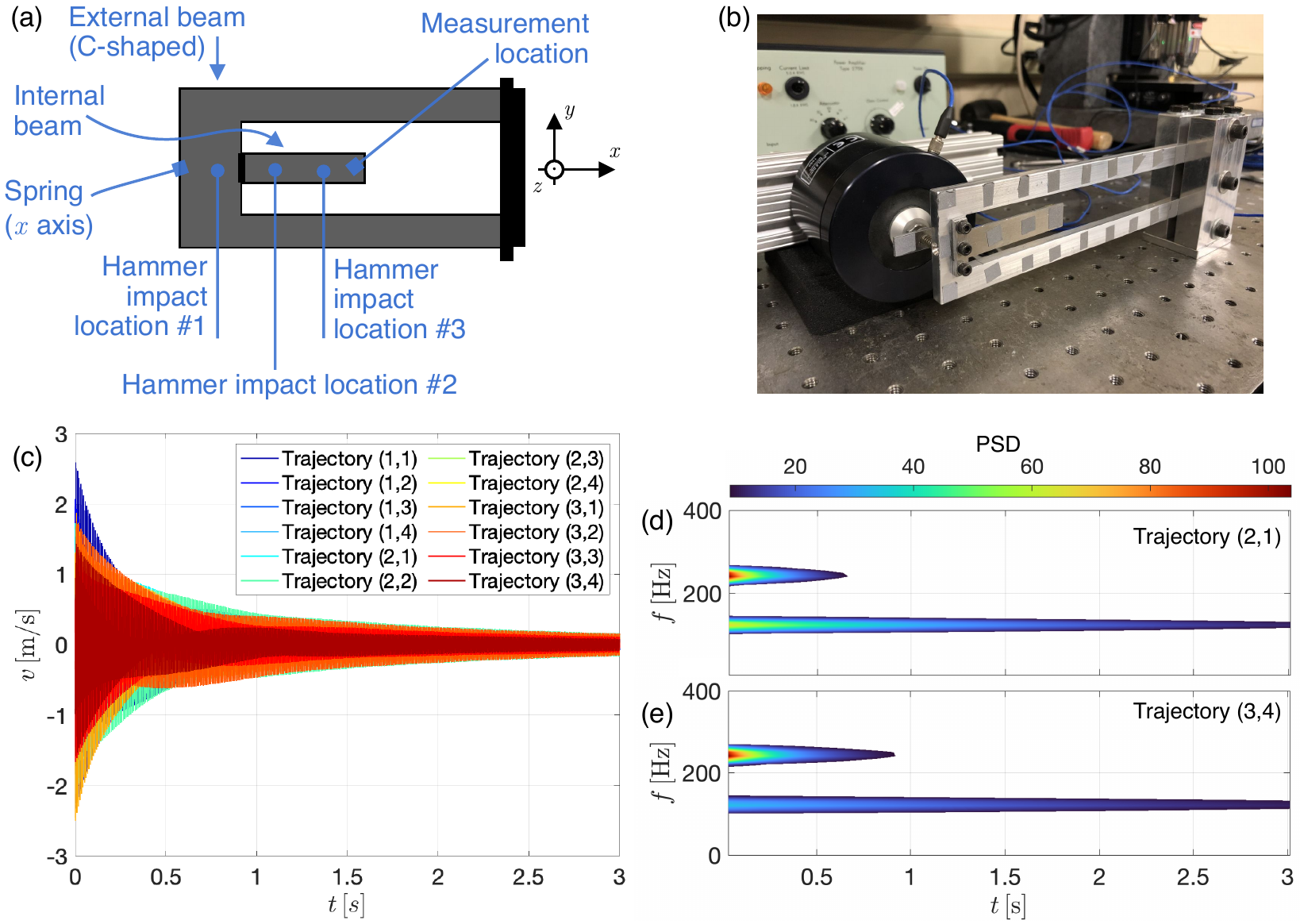}
    \caption{(a,b) Schematics and lab photo for the resonant tester. (c) Velocity time series of the inner beam tip. (d,e) Power spectral density computed via short-time Fourier transform (spectrogram) for two decaying responses of hammer impact tests.}
    \label{fig:S3_IR1}
\end{figure}

We validate our displacement-based SSM-reduced model on the data from the accelerometer located at $77$ [mm] from the left end of the beam, as shown in Fig. \ref{fig:S3_BRB}(a). This signal is reconstructed from the model by interpolating in the grid to obtain the accelerometer location and differentiating in time. The NMTE error is only $1.3\%$ and the model practically denoises the signal, as seen in Fig. \ref{fig:S3_BRB}(f,g). A further validation in Fig. \ref{fig:S3_BRB}(h,i) compares the instantaneous decay properties of the data-driven model (\ref{eq:RDBRB}) to those extracted using the Peak Finding and Fitting \cite{Jin2020,TRC2019_PartII} signal processing technique. There is close agreement among these curves, especially in the strongly nonlinear oscillation regime.

\subsection{Impacts on an internally resonant tester structure}
\label{sec:S3_IR}
Our final example is the resonant tester shown in Fig. \ref{fig:S3_IR1}(a,b). It consists of two beam-like parts made of aluminum 6061-T6, where the external beam is C-shaped and clamped to the ground on one side, while the internal beam is jointed to the external one via bolts three bolts: two side bolts are torqued to $1.36$ [Nm] for structural integrity whereas the middle bolt is torqued to $0.45$ [Nm] for enhanced frictional slip and associated nonlinearities. Additionally, a linear spring (Model $\#$1NCH2, Grainger, Inc.) connects the tip of the external beam to a fixed rigid frame in the direction of the z-axis. The system possesses an internal $1:2$ resonance between its slowest transverse bending modes, whose frequencies indeed clock at $122.4$ [Hz] and $243.4$ [Hz]. We consider transverse vibrations in the out of plane direction---the z-axis in Fig. \ref{fig:S3_IR1}(a). The available observable is the velocity of the inner beam tip, measured via laser scanner vibrometry (PSV400, Polytec Inc.). Transient vibrations are recorded for 3 seconds at a sampling rate of $5120$ [Hz].

A modally-tuned impulse hammer (PCB 086C01, PCB Piezotronics, Inc.) is used to excite transverse vibrations from three different impact locations in Fig. \ref{fig:S3_IR1}(c), so that our dataset features $12$ trajectories (four per impact location), shown in Fig. \ref{fig:S3_IR1}(c). We label these trajectories as $(j,l)$ where $j$ refers to the location and $l$ to the test number. Time-frequency analyses of the velocity signals, two of which are reported in Fig. \ref{fig:S3_IR1}(d,e), show that only the two slowest frequencies are present in the signal, so that the time responses can be well approximated by the slowest four-dimensional SSM of the system. The impact locations, the hammer tip and the forcing amounts were selected to achieve a sufficient trajectory diversity in the dataset without exciting further structural modes. For constructing an SSM-reduced model, we truncate the velocity signals after the hammer impact, use ten trajectories for training and leave two trajectories for testing.
\begin{figure}[t!]
    \centering
    \includegraphics[width=1\textwidth]{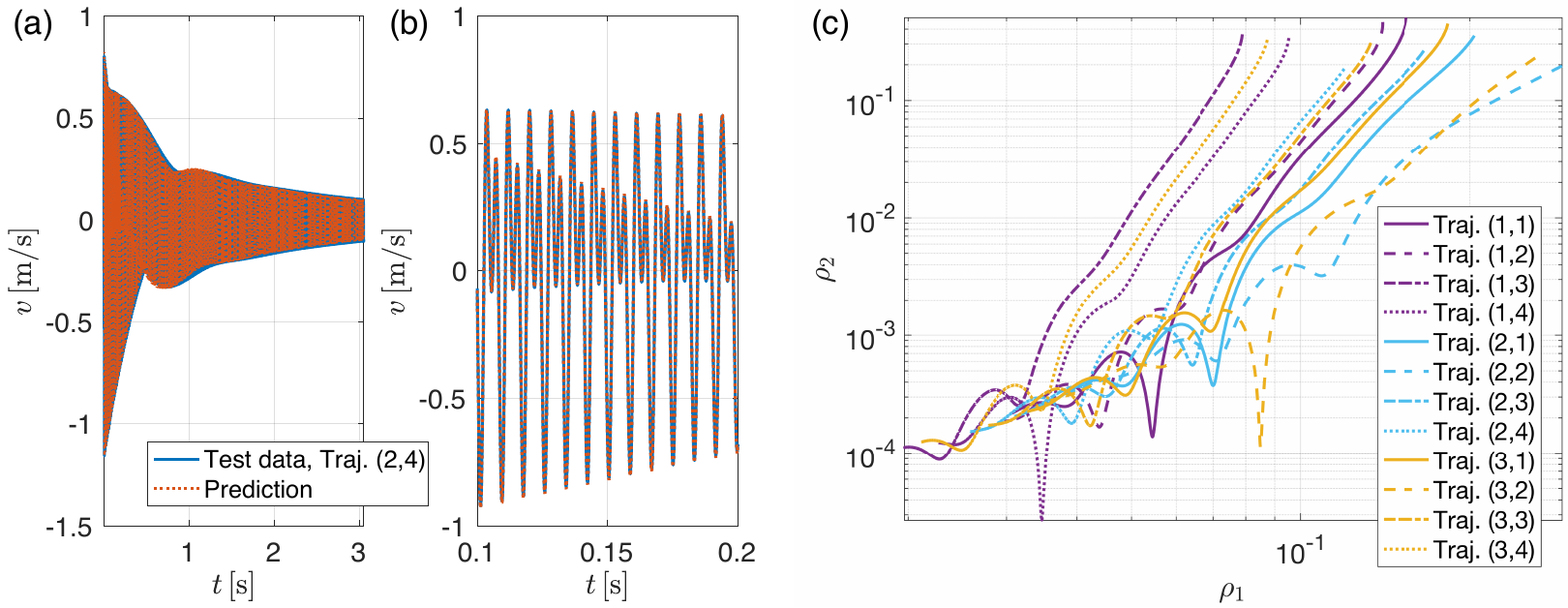}
    \caption{(a,b) Data-driven reduced-order model performances in reconstructing a testing trajectory. (c) Normal form amplitudes decays for the slow $\rho_1$ and fast $\rho_2$ modes, for all available trajectories in the dataset.}
    \label{fig:S3_IR2}
\end{figure}

The minimal embedding dimensions (nine for a four-dimensional manifold) fails to produce accurate reduced-order models (the NMTE error amounts to more than $8\%$). We therefore augment the delay embedding space so that each embedding vector captures approximately two cycles of the slowest oscillation. This procedure yields a $94$-dimensional delay embedding space. The result of our identification remains robust if we consider more embedding dimensions. A linear approximation to the embedded SSM has a good accuracy and our automated normal form algorithm, after estimating linearized eigenvalues, identifies a resonance among them. Defining $\psi = \theta_2-\theta_1$, we obtain from \texttt{SSMLearn} the cubic SSM-reduced polar normal form
\begin{equation}\label{eq:IRnf}\arraycolsep=1.4pt
\begin{array}{rl}\dot{\rho}_{1}&= -0.4228\rho_{1}-19.94 \rho^{3}_{1}+3.514 \rho_{1}\rho^{2}_{2}+\mathrm{Re}((0.08706-0.2427i) \rho_{2}\rho_{1} e^{i\psi})\\ & = - \alpha_1(\rho_1,\rho_2,\psi)\rho_1 \\
\dot{\rho}_{2} & = -3.155\rho_{2}-18.91\rho^{2}_{1}\rho_{2}-15.08 \rho^{3}_{2}+\mathrm{Re}((1.726-0.3342i) \rho^{2}_{1} e^{-i\psi})\\ & = - \alpha_2(\rho_1,\rho_2,\psi)\rho_2, \\
\rho_{1}\dot{\theta}_{1} &=+769.0\rho_{1}-59.56 \rho^{3}_{1}-0.5460 \rho^{2}_{2}\rho_{1}+\mathrm{Im}((0.08706-0.2427i) \rho_{2}\rho_{1} e^{i\psi}) \\ & = \omega_1(\rho_1,\rho_2,\psi)\rho_1,\\ 
\rho_{2}\dot{\theta}_{2} & =+1529\rho_{2}-31.26 \rho^{2}_{1}\rho_{2}-28.65 \rho^{3}_{2}+\mathrm{Im}((1.726-0.3342i) \rho^{2}_{1} e^{-i\psi})\\ & = \omega_2(\rho_1,\rho_2,\psi)\rho_2.
\end{array}
\end{equation}
\begin{figure}[t!]
    \centering
    \includegraphics[width=1\textwidth]{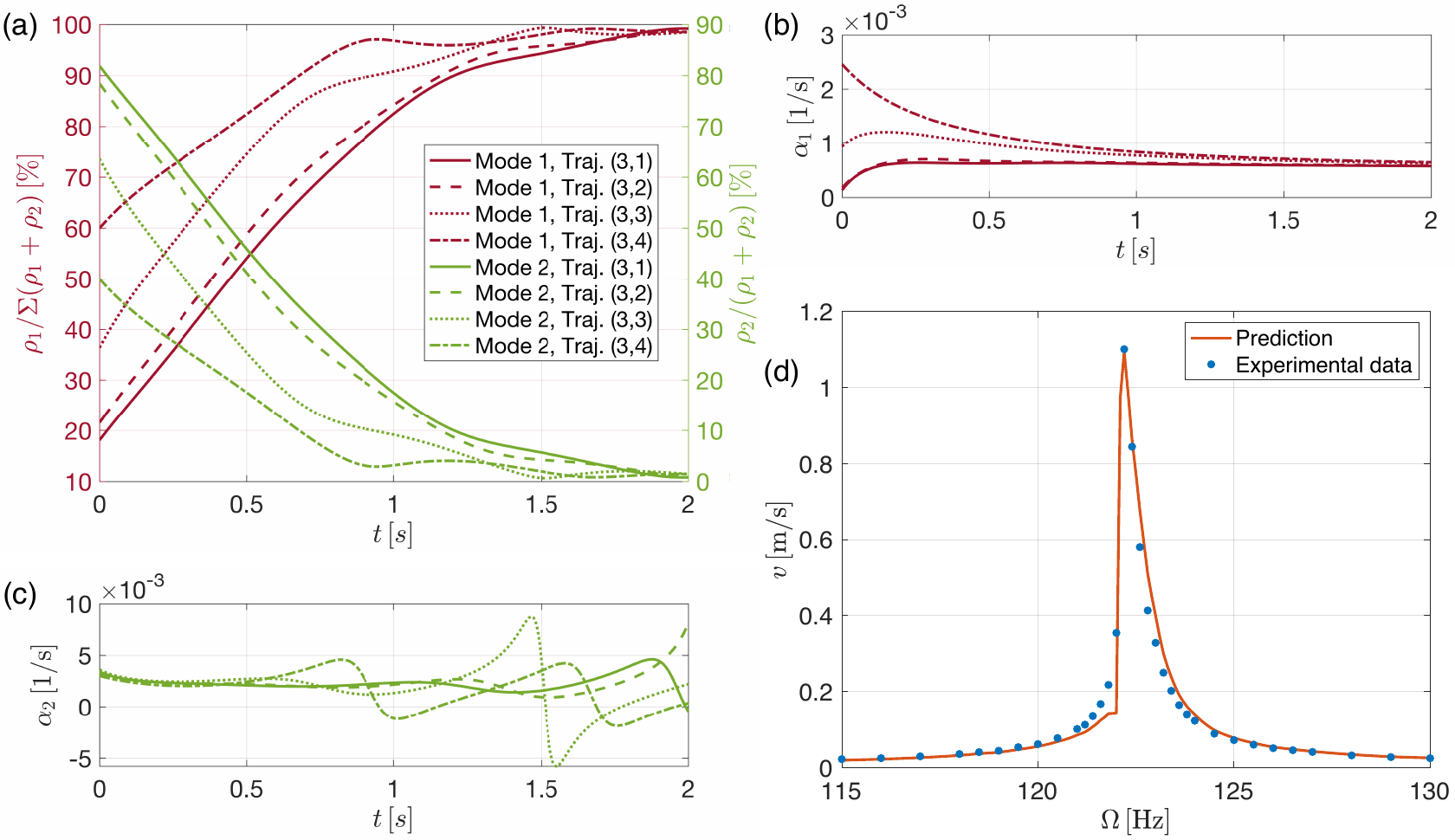}
    \caption{(a) Energy repartition in the resonant tester following a hammer impact on the third location. (b,c) The trend of instantaneous (or effective) damping of the normal form dynamics on the first two seconds of decays related to the third impact location. (d) Forced frequency response from experimental measurements and from analytical predictions based on the SSM-reduced model constructed by \texttt{SSMLearn}.}
    \label{fig:S3_IR3}
\end{figure}
This data-driven model reconstructs both test trajectories with an average $1.2$ $\%$ NMTE error, cf. Fig. \ref{fig:S3_IR2}(a,b). The decay of the slow mode amplitude $\rho_1$ and that of the fast one $\rho_2$ are shown in Fig. \ref{fig:S3_IR2}(c). Due to modal interactions, these decays are not monotone. From the plot, we notice a great diversity of decays depending on the impact location, and location three (the closest to the inner beam tip) is characterized by the highest amplitudes variability. 

Figure \ref{fig:S3_IR3}(a) shows energy repartition among the modes for the third impact location. This repartition is defined as the instantaneous ratio between the amplitude of a mode and the amplitude sum. Clearly, the slow mode tends to accumulate energy over time, while the fast mode dissipates it quickly. These trends are not monotonic, showing simultaneous and opposite changes in growth/decay rates, which implies that the faster mode is absorbing energy from the slower one. This behavior can also be inferred by the instantaneous properties illustrated in \ref{fig:S3_IR3}(b,c). The uncoupled limit of the oscillators suggests that the modes admit frequency softening and damping intensification when the oscillation amplitude increase. This is consistent with typical observations of jointed structures \cite{Brake2018}. In particular, the fast mode is coupled to the slow one and its damping undergoes consistent variation, becoming also negative for some times \cite{Sapsis2012}. Note that nonlinearity and coupling can be reduced at higher bolt torques, which, in the beam assembly used here, corresponds to tightening of the middle bolt. Coupling revealed by \texttt{SSMLearn} suggests that nearly decoupled modal oscillator models employed elsewhere \cite{Eriten2013,Segalman2015} are only valid for high bolt torques and small frictional slip, i.e., weak contact nonlinearities.

In addition to measuring decaying vibrations, we also perform some forced testing near the linearized frequency of the slow mode. We trigger forced responses in near-resonance with the slow (first bending) mode by using the Br\"uel $\&$ Kj{\ae}r 4810 shaker shown in Fig. \ref{fig:S3_IR1}(b), mounted on one end of the linear spring, and acquire velocity response from the tip of the inner beam by laser vibrometry. We also monitor the amplitude of shaker tip velocities and keep them constant while sweeping the frequencies around the first bending mode. In that sense, the response we obtain can be seen as transmissibilty rather than a classic frequency response curve, with forcing amplitudes kept constant throughout frequency sweeps. Starting from forced velocity time histories, we estimate the normal form forcing to be added to the vector field (\ref{eq:IRnf}) as in (\ref{eq:SSMdynf}). The resulting predictions are in very good agreement with experimental measurements in this weakly nonlinear regime, as shown in Fig. \ref{fig:S3_IR3}(d). Deeper analyses on forced responses are currently under investigation. Thanks to feedback loops used to track forcing, forced response curves can be extracted with improved accuracy, especially at non-linearizable amplitudes.

\section{Conclusion}
We have reviewed a general methodology for constructing sparse reduced-order models for potentially high-dimensional, nonlinear mechanical systems from data. Our approach constructs normal forms on attracting spectral submanifolds (SSMs), which are the smoothest nonlinear continuation of spectral subspaces of the linearized dynamics. Implemented in the publicly available \textsc{Matlab}\textsuperscript{\textregistered} code \texttt{SSMLearn} \cite{SSMLearn}, our algorithm takes generic observable data as input to identify robust and predictive nonlinear models that also capture for non-linearizable dynamics. SSM theory offers a systematic basis for model reduction and allows a simplification of the reduced dynamics via normal forms, which are particularly insightful for mechanical systems. Indeed, SSM-reduced models can handle multi-modal interactions, identify amplitude-dependent damping and frequency, and predict the forced structural response.

We have illustrated SSM-reduced modeling in numerical and experimental case studies, featuring different types of observables, nonlinearities and SSM dimensions. Specifically, we have discussed different dynamical regimes and the relevance of slow SSMs in a chain of oscillators, derived a reduced-order model from digital image correlation measurements of the Brake-Reuss beam, and unfolded the internally resonant dynamics of a tester structure, also predicting forced responses. These examples were analyzed using the open-source \textsc{Matlab}\textsuperscript{\textregistered} package \texttt{SSMLearn} \cite{SSMLearn} that performs data-driven, SSM-based model reduction starting from vibrations data. This algorithm only requires a minimal number of input parameters: the SSM dimension, the polynomial order for SSM parametrization and the polynomial order of the reduced dynamics. The SSM dimension is either known a priori from targeted experiments (e.g. resonance decay) or can be estimated via time-frequency signal processing analysis of the input data. This makes our method a parametric complement to non-parametric identification tools. Polynomial orders can be adjusted to improve accuracy, noting that excessive orders may lead to overfitting. With the help of the numerical continuation core \textsc{coco} \cite{Dankowicz2013} included in \texttt{SSMLearn} \cite{SSMLearn}, users can compute forced response curves or design nonlinear control strategies from the identified nonlinear models. 

Further examples, both numerical and experimental, with detailed code are available in the \texttt{SSMLearn} repository \cite{SSMLearn}, which is also suitable for high-dimensional fluid flows and fluid structure interaction problems \cite{Cenedese2021}. Current limitations of the present approach include weaker performance for large forcing amplitudes. These appear, for example, in the Brake-Reuss beam experiments of \cite{TRC2019_PartI,TRC2019_PartII}, which we expect to capture only with a more refined forced-reduced dynamics and improved calibration procedures. The same requirement holds for more complicated forcing types (e.g., quasi-periodic or random), which are relevant in structural dynamics. We are addressing these challenges in ongoing work.

%%%%%%%%%% Insert bibliography here %%%%%%%%%%%%%%


\begin{thebibliography}{10}

\bibitem{Peeters2011a}
M.~Peeters, G.~Kerschen, and J.C. Golinval.
 Dynamic testing of nonlinear vibrating structures using nonlinear
  normal modes.
 {\em Journal of Sound and Vibration}, 330(3):486--509, 2011.

\bibitem{Peeters2011b}
M.~Peeters, G.~Kerschen, and J.C. Golinval.
 Modal testing of nonlinear vibrating structures based on nonlinear
  normal modes: experimental demonstration.
 {\em Mechanical Systems and Signal Processing}, 25(4):1227--1247,
  2011.

\bibitem{Kurt2014}
M.~Kurt, M.~Eriten, D.M. McFarlan, L.A. Bergman, and A.F. Vakakis.
 Strongly nonlinear beats in the dynamics of an elastic system with a
  strong local stiffness nonlinearity: Analysis and identification.
 {\em Journal of Sound and Vibration}, 333(7):2054--2072, 2014.

\bibitem{Ehrhardt2016}
D.A. Ehrhardt and M.S. Allen.
 Measurement of nonlinear normal modes using multi-harmonic stepped
  force appropriation and free decay.
 {\em Mechanical Systems and Signal Processing}, 76-77:612--633, 2016.

\bibitem{Renson2016b}
L.~Renson, A.~Gonzalez-Buelga, D.A.W. Barton, and S.A. Neild.
 Robust identification of backbone curves using control-based
  continuation.
 {\em Journal of Sound and Vibration}, 367:145--158, 2016.

\bibitem{Moore2019}
K.J. Moore, M.~Kurt, M.~Eriten, D.M. McFarlan, L.A. Bergman, and A.F. Vakakis.
 Direct detection of nonlinear modal interactions from time series
  measurements.
 {\em Mechanical Systems and Signal Processing}, 125:311--329, 2019.
 Exploring nonlinear benefits in engineering.

\bibitem{TRC2019_PartI}
W.~Chen, D.~Jana, A.~Singh, M.~Jin, M.~Cenedese, G.~Kosova, M.W.R. Brake, C.W.
  Schwingshackl, S.~Nagarajaiah, K.~Moore, and J.P. Noël.
 Measurement and identification of the nonlinear dynamics of a jointed
  structure using full-field data; {P}art {I} - {M}easurement of nonlinear
  dynamics.
 {\em Mechanical Systems and Signal Processing}, 166:108401, 2022.

\bibitem{TRC2019_PartII}
M.~Jin, G.~Kosova, M.~Cenedese, W.~Chen, D.~Jana, A.~Singh, M.W.R. Brake, C.W.
  Schwingshackl, S.~Nagarajaiah, K.~Moore, and J.P. Noël.
 Measurement and identification of the nonlinear dynamics of a jointed
  structure using full-field data; {P}art {II} - {N}onlinear system
  identification.
 {\em Mechanical Systems and Signal Processing}, 166:108402, 2022.

\bibitem{Awrejcewicz2004}
J.~Awrejcewicz, V.A. Krys'ko, and A.F. Vakakis.
 {\em Order Reduction by Proper Orthogonal Decomposition (POD)
  Analysis}, pages 279--320.
 Springer, Berlin, Heidelberg, 2004.

\bibitem{Holmes2012}
P.~Holmes, J.L. Lumley, G.~Berkooz, and C.W. Rowley.
 {\em Turbulence, Coherent Structures, Dynamical Systems and
  Symmetry}.
 Cambridge Monographs on Mechanics. Cambridge University Press, 2
  edition, 2012.

\bibitem{Lu2019}
K.~Lu, Y.~Jin, Y.~Chen, Y.~Yang, L.~Hou, Z.~Zhang, Z.~Li, and C.~Fu.
 Review for order reduction based on proper orthogonal decomposition
  and outlooks of applications in mechanical systems.
 {\em Mechanical Systems and Signal Processing}, 123:264--297, 2019.

\bibitem{Hijazi2020}
S.~Hijazi, G.~Stabile, A.~Mola, and G.~Rozza.
 Data-driven {POD}-{G}alerkin reduced order model for turbulent flows.
 {\em Journal of Computational Physics}, 416:109513, 2020.

\bibitem{Schmid2010}
P.J. Schmid.
 Dynamic mode decomposition of numerical and experimental data.
 {\em Journal of Fluid Mechanics}, 656:5–28, 2010.

\bibitem{Kutz2016}
J.N. Kutz, S.L. Brunton, B.W. Brunton, and J.L. Proctor.
 {\em Dynamic Mode Decomposition}.
 Society for Industrial and Applied Mathematics, Philadelphia, PA,
  2016.

\bibitem{Chen2012}
K.K. Chen, J.H. Tu, and C.W. Rowley.
 Variants of dynamic mode decomposition: Boundary condition,
  {K}oopman, and {F}ourier analyses.
 {\em Annual Review of Fluid Mechanics}, 22(6):887--915, 2012.

\bibitem{Williams2015}
M.O. Williams, C.W. Rowley, and I.G. Kevrekidis.
 A kernel-based method for data-driven koopman spectral analysis.
 {\em Journal of Computational Dynamics}, 2(2):247--265, 2015.

\bibitem{Rowley2017}
C.W. Rowley and S.T.M. Dawson.
 Model reduction for flow analysis and control.
 {\em Annual Review of Fluid Mechanics}, 49(1):387--417, 2017.

\bibitem{Alla2017}
A.~Alla and J.N. Kutz.
 Nonlinear model order reduction via dynamic mode decomposition.
 {\em SIAM Journal on Scientific Computing}, 39(5):B778--B796, 2017.

\bibitem{Lusch2018}
B.~Lusch, J.N. Kutz, and S.L. Brunton.
 Deep learning for universal linear embeddings of nonlinear dynamics.
 {\em Nature Communications}, 9:4950, 2018.

\bibitem{Budisic2012}
M.~Budišić, R.~Mohr, and I.~Mezić.
 Applied koopmanism.
 {\em Chaos: An Interdisciplinary Journal of Nonlinear Science},
  22(4):047510, 2012.

\bibitem{Mezic2013}
I.~Mezić.
 Analysis of fluid flows via spectral properties of the koopman
  operator.
 {\em Annual Review of Fluid Mechanics}, 45(1):357--378, 2013.

\bibitem{Bagheri2013}
S.~Bagheri.
 Koopman-mode decomposition of the cylinder wake.
 {\em Journal of Fluid Mechanics}, 726:596--623, 2013.

\bibitem{Brunton2016b}
S.L. Brunton, B.W. Brunton, J.L. Proctor, and J.N. Kutz.
 Koopman invariant subspaces and finite linear representations of
  nonlinear dynamical systems for control.
 {\em PLOS ONE}, 11:1--19, 2016.

\bibitem{Page2019}
J.~Page and R.R. Kerswell.
 Koopman mode expansions between simple invariant solutions.
 {\em Journal of Fluid Mechanics}, 879:1–27, 2019.

\bibitem{Kaiser2021}
E.~Kaiser, J.N. Kutz, and S.L. Brunton.
 Data-driven discovery of koopman eigenfunctions for control.
 {\em Machine Learning: Science and Technology}, 2(3):035023, 2021.

\bibitem{Page2018}
J.~Page and R.R. Kerswell.
 Koopman analysis of burgers equation.
 {\em Phys. Rev. Fluids}, 3:071901, Jul 2018.

\bibitem{Pearson1901}
K.F.R.S. Pearson.
 {LIII}. {O}n lines and planes of closest fit to systems of points in
  space.
 {\em The London, Edinburgh, and Dublin Philosophical Magazine and
  Journal of Science}, 2(11):559--572, 1901.

\bibitem{Bishop2006}
C.M. Bishop.
 {\em Pattern Recognition and Machine Learning}.
 Information Science and Statistics. Springer-Verlag New York, 2006.

\bibitem{Roweis200}
S.T. Roweis and L.K. Saul.
 Nonlinear dimensionality reduction by locally linear embedding.
 {\em Science}, 290(5500):2323--2326, 2000.

\bibitem{Bengio2013}
Y.~{Bengio}, A.~{Courville}, and P.~{Vincent}.
 Representation learning: a review and new perspectives.
 {\em IEEE Transactions on Pattern Analysis and Machine Intelligence},
  35(8):1798--1828, 2013.

\bibitem{Loiseau2020}
J.-C. Loiseau, S.L. Brunton, and B.R. Noack.
 From the {POD}-{G}alerkin method to sparse manifold models.
 In P.~Benner, editor, {\em Handbook of Model Order Reduction, Volume
  {II}: Applications}, pages 1--47. De Gruyter GmbH, 2019.

\bibitem{Bengio2016}
I.~Goodfellow, Y.~Bengio, and A.~Courville.
 {\em Deep Learning}.
 MIT Press, 2016.

\bibitem{Champion2019}
K.~Champion, B.~Lusch, J.N. Kutz, and S.L. Brunton.
 Data-driven discovery of coordinates and governing equations.
 {\em Proceedings of the National Academy of Sciences},
  116(45):22445--22451, 2019.

\bibitem{Tibshirani1996}
R.~Tibshirani.
 Regression shrinkage and selection via the lasso.
 {\em Journal of the Royal Statistical Society. Series B
  (Methodological)}, 58(1):267--288, 1996.

\bibitem{Brunton2016}
S.L. Brunton, J.L. Proctor, and J.N. Kutz.
 Discovering governing equations from data by sparse identification of
  nonlinear dynamical systems.
 {\em Proceedings of the National Academy of Sciences},
  113(15):3932--3937, 2016.

\bibitem{Abdessalem2019}
A.~{Ben Abdessalem}, N.~Dervilis, D.J. Wagg, and K.~Worden.
 Model selection and parameter estimation of dynamical systems using a
  novel variant of approximate bayesian computation.
 {\em Mechanical Systems and Signal Processing}, 122:364--386, 2019.

\bibitem{Billings2013}
S.A. Billings.
 {\em Nonlinear System Identification: NARMAX Methods in the Time,
  Frequency, and Spatio-Temporal Domains}.
 Springer, 2013.

\bibitem{Hartman2017}
D.~Hartman and L.K. Mestha.
 A deep learning framework for model reduction of dynamical systems.
 In {\em 2017 IEEE Conference on Control Technology and Applications
  (CCTA)}, pages 1917--1922, 2017.

\bibitem{Lui2019}
H.F.S. Lui and W.R. Wolf.
 Construction of reduced-order models for fluid flows using deep
  feedforward neural networks.
 {\em Journal of Fluid Mechanics}, 872:963–994, 2019.

\bibitem{Raissi2019}
M.~Raissi, P.~Perdikaris, and G.E. Karniadakis.
 Physics-informed neural networks: A deep learning framework for
  solving forward and inverse problems involving nonlinear partial differential
  equations.
 {\em Journal of Computational Physics}, 378:686--707, 2019.

\bibitem{Karniadakis2021}
G.E. Karniadakis, I.G. Kevrekidis, L.~Lu, P.~Perdikaris, S.~Wang, and L.~Yang.
 Physics-informed machine learning.
 {\em Nature Reviews Physics}, 3:422–440, 2021.

\bibitem{Haller2016}
G.~Haller and S.~Ponsioen.
 Nonlinear normal modes and spectral submanifolds: existence,
  uniqueness and use in model reduction.
 {\em Nonlinear Dynamics}, 86(3):1493--1534, 2016.

\bibitem{Haller2017}
G.~Haller and S.~Ponsioen.
 Exact model reduction by a slow-fast decomposition of nonlinear
  mechanical systems.
 {\em Nonlinear Dynamics}, 90:617--647, 2017.

\bibitem{Szalai2017}
R.~Szalai, D.~Ehrhardt, and G.~Haller.
 Nonlinear model identification and spectral submanifolds for
  multi-degree-of-freedom mechanical vibrations.
 {\em Proceedings of the Royal Society of London A: Mathematical,
  Physical and Engineering Sciences}, 473:20160759, 2017.

\bibitem{Jain2018}
S.~Jain, P.~Tiso, and G.~Haller.
 Exact nonlinear model reduction for a von {K}ármán beam: slow-fast
  decomposition and spectral submanifolds.
 {\em Journal of Sound and Vibration}, 423:195--211, 2018.

\bibitem{Ponsioen2018}
S.~Ponsioen, T.~Pedergnana, and G.~Haller.
 Automated computation of autonomous spectral submanifolds for
  nonlinear modal analysis.
 {\em Journal of Sound and Vibration}, 420:269--295, 2018.

\bibitem{Breunung2018}
T.~Breunung and G.~Haller.
 Explicit backbone curves from spectral submanifolds of forced-damped
  nonlinear mechanical systems.
 {\em Proceedings of the Royal Society of London A: Mathematical,
  Physical and Engineering Sciences}, 474:20180083, 2018.

\bibitem{Ponsioen2019}
S.~Ponsioen, T.~Pedergnana, and G.~Haller.
 Analytic prediction of isolated forced response curves from spectral
  submanifolds.
 {\em Nonlinear Dynamics}, 98:2755--2773, 2019.

\bibitem{Ponsioen2020}
S.~Ponsioen, S.~Jain, and G.~Haller.
 Model reduction to spectral submanifolds and forced-response
  calculation in high-dimensional mechanical systems.
 {\em Journal of Sound and Vibration}, 488:115640, 2020.

\bibitem{Jain2021}
S.~Jain and G.~Haller.
 How to compute invariant manifolds and their reduced dynamics in
  high-dimensional finite-element models?
 {\em Nonlinear Dynamics}, in press, 2021.

\bibitem{SSMTool2021}
S.~Jain, T.~Thurnher, M.~Li, and G.~Haller.
\texttt{SSMTool-2.0}: Computation of invariant manifolds \& their reduced
  dynamics in high-dimensional mechanics problems,  \href{https://doi.org/10.5281/zenodo.4614202}{www.georgehaller.com}, 2021.

\bibitem{Szalai2020}
R.~Szalai.
Invariant spectral foliations with applications to model order
  reduction and synthesis.
{\em Nonlinear Dynamics}, 101:2645--2669, 2020.

\bibitem{Cenedese2021}
M.~Cenedese, J.~Ax\r{a}s, B.~B\"auerlein, K.~Avila, and G.~Haller.
 Data-driven modeling and prediction of non-linearizable dynamics via
  spectral submanifolds.
 Submitted, 2021.

\bibitem{SSMLearn}
M.~Cenedese, J.~Ax\r{a}s, and G.~Haller.
\texttt{SSMLearn}: Data-driven reduced order models for nonlinear dynamical systems. To be released, 2021.

\bibitem{Noel2017}
J.P. Noël and G.~Kerschen.
 Nonlinear system identification in structural dynamics: 10 more years
  of progress.
 {\em Mechanical Systems and Signal Processing}, 83:2--35, 2017.

\bibitem{Feldman2011}
M.~Feldman.
 Hilbert transform in vibration analysis.
 {\em Mechanical Systems and Signal Processing}, 25(3):735--802, 2011.

\bibitem{Jin2020}
M.~Jin, W.~Chen, M.R.W. Brake, and H.~Song.
 Identification of instantaneous frequency and damping from transient
  decay data.
 {\em Journal of Vibration and Acoustics}, 142(5):051111, 2020.

\bibitem{Moore2018}
K.J. Moore, M.~Kurt, M.~Eriten, D.M. McFarland, L.A. Bergman, and A.F. Vakakis.
 Wavelet-bounded empirical mode decomposition for measured time series
  analysis.
 {\em Mechanical Systems and Signal Processing}, 99:14--29, 2018.

\bibitem{Cabre2003a}
X.~Cabré, E.~Fontich, and R.~De~la lLave.
 The parameterization method for invariant manifolds {I}: manifolds
  associated to non-resonant subspaces.
 {\em Indiana University Mathematics Journal}, 52(2):283--328, 2003.

\bibitem{Cabre2003b}
X.~Cabré, E.~Fontich, and R.~De~la lLave.
 The parameterization method for invariant manifolds {II}: regularity
  with respect to parameters.
 {\em Indiana University Mathematics Journal}, 52(2):329--360, 2003.

\bibitem{Cabre2005}
X.~Cabré, E.~Fontich, and R.~De~la lLave.
 The parameterization method for invariant manifolds {III}: overview
  and applications.
 {\em Journal of Differential Equations}, 218(2):444--515, 2005.

\bibitem{Haro2016}
A.~Haro, M.~Canadell, J.~Figueras, and J.M. Luque, A.~Mondelo.
 {\em The Parameterization Method for Invariant Manifolds}, volume 195
  of {\em Applied Mathematical Sciences}.
 Springer International Publishing, 1st edition, 2016.

\bibitem{Fenichel1974}
N.~Fenichel.
 Asymptotic stability with rate conditions.
 {\em Indiana University Mathematics Journal}, 23(12):1109--1137,
  1974.

\bibitem{Sauer1991}
T.~Sauer, J.A. Yorke, and M.~Casdagli.
 Embedology.
 {\em Journal of Statistical Physics}, 65:579–616, 1997.

\bibitem{Takens1981}
F.~Takens.
 Detecting strange attractors in turbulence.
 In D.~Rand and L.~Young, editors, {\em Dynamical Systems and
  Turbulence, Warwick 1980}, pages 366--381. Springer Berlin Heidelberg,
  Berlin, Heidelberg, 1981.

\bibitem{GH1983}
J.~Guckenheimer and P.J. Holmes.
 {\em Nonlinear Oscillations, Dynamical Systems, and Bifurcations of
  Vector Fields}, volume~42 of {\em Applied Mathematical Sciences}.
 Springer-Verlag New York, 1983.

\bibitem{Murdock2003}
J.~Murdock.
 {\em Normal Forms and Unfoldings for Local Dynamical Systems}.
 Springer Monographs in Mathematics. Springer-Verlag New York, 2003.

\bibitem{Poincare1892}
H.~Poincaré.
 {\em Les Méthodes Nouvelles de la Mécanique Céleste}.
 Gauthier-Villars et Fils, Paris, 1892.

\bibitem{Marsden1976}
J.E. Marsden and M.~McCracken.
 {\em The Hopf Bifurcation and Its Applications}, volume~19 of {\em
  Applied Mathematical Sciences}.
 Springer-Verlag New York, 1976.

\bibitem{Cenedese2019}
M.~Cenedese and G.~Haller.
 How do conservative backbone curves perturb into forced responses?
  {A} {M}elnikov function analysis.
 {\em Proceedings of the Royal Society of London A: Mathematical,
  Physical and Engineering Sciences}, 476:20190494, 2020.

\bibitem{Dankowicz2013}
H.~Dankowicz and F.~Schilder.
 {\em Recipes for Continuation}.
 Society for Industrial and Applied Mathematics, 2013.

\bibitem{Brake2018}
M.R.W. Brake.
 {\em The mechanics of jointed structures: recent research and open
  challenges for developing predictive models for structural dynamics}.
 Springer International Publishing, 2018.

\bibitem{Sapsis2012}
T.P. Sapsis, D.D. Quinn, A.F. Vakakis, and L.A. Bergman.
 Effective stiffening and damping enhancement of structures with
  strongly nonlinear local attachments.
 {\em Journal of Vibration and Acoustics}, 134(1):011016, 2012.

\bibitem{Eriten2013}
M.~Eriten, M.~Kurt, G.~Luo, D.M. D.~McFarland, L.A. Bergman, and A.F. Vakakis.
 Nonlinear system identification of frictional effects in a beam with
  a bolted joint connection.
 {\em Mechanical Systems and Signal Processing}, 39(1):245--264, 2013.

\bibitem{Segalman2015}
D.J. Segalman, M.S. Allen, M. Eriten, and K. Hoppman. 
Experimental assessment of joint-like modal models for structures.
{\em Proceedings of the ASME 2015 International Design Engineering Technical Conferences and Computers and Information in Engineering Conference. Volume 8: 27th Conference on Mechanical Vibration and Noise}. Boston, Massachusetts, USA. August 2–5, 2015. V008T13A025. ASME.

\end{thebibliography}
\end{document}